\documentclass[12pt]{article}
\newtheorem{la}{Lemma}
\newtheorem{thm}[la]{Theorem}
\newtheorem{cor}[la]{Corollary}
\def\qed{\hspace*{\fill}\vrule height6pt width4pt depth0pt\medskip}
\def\quax{\hskip 1cm}

\long\def\longdelete#1{}

\def\ip{{\rm ip}}
\def\nc{{\rm nc}}

\begin{document}
\title{{\scshape\normalsize Mathematics Division, National Center
                            for Theoretical Sciences at Taipei}\\
{\scshape\large NCTS/TPE-Math Technical Report 2004-004}\\~\\
{\bf Isometric-path numbers of block graphs}}
\author{Jun-Jie Pan\thanks{Department of Applied Mathematics,
              National Chiao Tung University, Hsinchu 30050, Taiwan.}
        \and
        Gerard J. Chang\thanks{Department of Mathematics,
              National Taiwan University, Taipei 10617, Taiwan.
        E-mail: gjchang @math.ntu.edu.tw.
        Supported in part by the National Science Council
              under grant NSC90-2115-M002-024.
        Member of Mathematics Division, National Center for
              Theoretical Sciences at Taipei.}
       }
\date{Revised on March 10, 2004} 
\maketitle

\begin{abstract}
An isometric path between two vertices in a graph $G$ is a
shortest path joining them. The isometric-path number of $G$,
denoted by $\ip(G)$, is the minimum number of isometric paths
required to cover all vertices of $G$. In this paper, we determine
exact values of isometric-path numbers of block graphs. We also
give a linear-time algorithm for finding the corresponding paths.

\bigskip

\bigskip

\end{abstract}

\baselineskip=24pt

\section{Introduction}

An {\it isometric path} between two vertices in a graph $G$ is a
shortest path joining them. The {\it isometric-path number} of
$G$, denoted by $\ip(G)$, is the minimum number of isometric paths
required to cover all vertices of $G$. This concept has a close
relationship with the game of cops and robbers described as
follows. The game is played by two players, the {\it cop} and the
{\it robber}, on a graph. The two players move alternatively,
starting with the cop. Each player's first move consists of
choosing a vertex at which to start. At each subsequent move, a
player may choose either to stay at the same vertex or to move to
an adjacent vertex. The object for the cop is to catch the robber,
and for the robber is to prevent this from happening. Nowakowski
and Winkler \cite{nw83} and Quilliot \cite{q78} independently
proved that the cop wins if and only if the graph can be reduced
to a single vertex by successively removing pitfalls, where a {\it
pitfall} is a vertex whose closed neighborhood is a subset of the
closed neighborhood of another vertex. As not all graphs are
cop-win graphs, Aigner and Fromme \cite{af84} introduced the
concept of the {\it cop-number} of a general graph $G$, denoted by
$c(G)$, which is the minimum number of cops needed to put into the
graph in order to catch the robber . On the way to giving an upper
bound for the cop-numbers of planar graphs, they showed that a
single cop moving on an isometric path $P$ guarantee that after a
finite number of moves the robber will be immediately caught if he
moves onto $P$. Observing this fact, Fitzpatrick \cite{f97} then
introduced the concept of isometric-path cover and pointed out
that $c(G) \le \ip(G)$.

The isometric-path number of the Cartesian product $P_{n_1} \times
P_{n_2} \times \ldots \times P_{n_d}$ has been studied in the
literature. Fitzpatrick \cite{f99} gave bounds for the case when
$n_1=n_2=\ldots=n_d$. Fisher and Fitzpatrick \cite{ff01} gave
exact values for the case $d=2$. Fitzpatrick et al \cite{fnhc01}
gave a lower bound, which is in fact the exact value if $d+1$ is a
power of $2$, for the case when $n_1=n_2=\ldots=n_d=2$.

The purpose of this paper is to give exact values of
isometric-path numbers of block graphs. We also give a linear-time
algorithm to find the corresponding paths. For technical reasons,
we consider a slightly more general problem as follows. Suppose
every vertex $v$ in the graph $G$ is associated with a
non-negative integer $f(v)$. We call such function $f$ a {\it
vertex labeling} of $G$. An {\it $f$-isometric-path cover} of $G$
is a family $\cal C$ of isometric paths such that the following
conditions hold.

\begin{description}
\item
(C1) If $f(v)=0$, then $v$ is in an isometric path in $\cal C$.

\item (C2) If $f(v)\ge 1$, then $v$ is an end vertex of at least
$f(v)$ isometric paths in $\cal C$, while the counting is twice if
$v$ itself is a path in $\cal C$.
\end{description}

The $f$-{\it isometric-path number} of $G$, denoted by $\ip_f(G)$,
is the minimum cardinality of an $f$-isometric-path cover of $G$.
It is clear that when $f(v)=0$ for all vertices $v$ in $G$, we
have $\ip(G)=\ip_f(G)$. The attempt of is paper is to determine
the $f$-isometric-path number of a block graph. Recall that a {\it
block graph} is a graph in which every block is a complete graph.
A {\it cut-vertex} of a graph is a vertex whose removal results in
a graph with more components than the original graph. It is
well-known that in a block graph all internal vertices of an
isometric path are cut-vertices.

\section{Block graphs}

In this section, we determine the $f$-isometric-path numbers for
block graphs $G$. Without loss of generality, we may assume that
$G$ is connected.

First, a useful lemma.

\begin{la} \label{la:1} 
Suppose $x$ is a non-cut-vertex of a block graph $G$
with a vertex labeling $f$.
If vertex labeling $f'$ is the same as $f$ except that $f'(x)=\max\{1,f(x)\}$,
then $\ip_f(G)=\ip_{f'}(G)$.
\end{la}
{\bf Proof.} As any internal vertex of an isometric path in a
block graph is a cut-vertex but $x$ not a cut-vertex, $x$ must be
an end vertex of any isometric path. It follows that a collection
$\cal C$ is an $f$-isometric-path cover if and only if it is an
$f'$-isometric-path cover. The lemma then follows. \qed

So, now we may assume that $f(v) \ge 1$ for all non-cut-vertices $v$ of $G$,
and call such a vertex labeling {\it regular}.
Now, we have the following theorem for the inductive step.

\begin{thm} \label{thm:2} 
Suppose $G$ is a block graph with a regular labeling $f$,
and $x$ is a non-cut-vertex in a block $B$
with exactly one cut-vertex $y$ or with no cut-vertex
in which case let $y$ be any vertex of $B-\{x\}$.
When $f(x)=1$, let $G'=G-x$ with a regular vertex labeling $f'$
which is the same as $f$ except $f'(y)=f(y)+1$.
When $f(x) \ge 2$, let $G'=G$ with a regular vertex labeling $f'$
which is the same as $f$ except $f'(x)=f(x)-1$ and $f'(y)=f(y)+1$.
Then $\ip_f(G)=\ip_{f'}(G')$.
\end{thm}
{\bf Proof.} We first prove that $\ip_f(G)\ge \ip_{f'}(G')$.
Suppose $\cal C$ is an optimal $f$-isometric-path cover of $G$.
Choose a path $P$ in $\cal C$ having $x$ as an end vertex. We
consider four cases.

{\bf Case 1.1.} $P=x$ and $f(x)=1$ (i.e., $G'=G-x$).

In this case, ${\cal C}'=({\cal C}-\{P\})\cup \{y\}$ is an
$f'$-isometric-path cover of $G'$. Hence, $\ip_f(G) = |{\cal C}|
\ge |{\cal C}'| \ge \ip_{f'}(G')$.

{\bf Case 1.2.} $P=x$ and $f(x)\ge 2$ (i.e., $G'=G$).

In this case, ${\cal C}'=({\cal C}-\{P\})\cup \{xy\}$ is an
$f'$-isometric-path cover of $G'$. Hence, $\ip_f(G) = |{\cal C}|
\ge |{\cal C}'| \ge \ip_{f'}(G')$.

{\bf Case 1.3.} $P=xz$ for some vertex $z$ in $B-\{x,y\}$.

In this case, ${\cal C}'=({\cal C}-\{P\})\cup \{yz\}$ is an
$f'$-isometric-path cover of $G'$. Hence, $\ip_f(G) = |{\cal C}|
\ge |{\cal C}'| \ge \ip_{f'}(G')$.

{\bf Case 1.4.} $P=xyQ$, where $Q$ contains no vertices in $B$.

In this case, ${\cal C}'=({\cal C}-\{P\})\cup \{yQ\}$ is an
$f'$-isometric-path cover of $G'$. Hence, $\ip_f(G) = |{\cal C}|
\ge |{\cal C}'| \ge \ip_{f'}(G')$.

Next, we prove that $\ip_f(G) \le \ip_{f'}(G')$. Suppose ${\cal
C}'$ is an optimal $f'$-isometric-path cover of $G'$. Choose a
path $P'$ in ${\cal C}'$ having $y$ as an end vertex. We consider
three cases.

{\bf Case 2.1.} $P'=yx$.

In this case, $G'=G$ and ${\cal C}=({\cal C'}-\{P'\})\cup \{x\}$
is an $f$-isometric-path cover of $G$. Hence, $\ip_f(G) \le |{\cal
C}| \le |{\cal C}'| = \ip_{f'}(G')$.

{\bf Case 2.2.} $P'=yz$ for some $z$ in $B-\{x,y\}$.

In this case, ${\cal C}=({\cal C}'-\{P'\})\cup \{xz\}$ is an
$f$-isometric-path cover of $G$. Hence, $\ip_f(G) \le |{\cal C}|
\le |{\cal C}'| = \ip_{f'}(G')$.

{\bf Case 2.3.} $P'=yQ$, where $Q$ contains no vertex in $B$.

In this case, ${\cal C}=({\cal C}'-\{P'\})\cup \{xyQ\}$ is an
$f$-isometric-path cover of $G$. Hence, $\ip_f(G) \le |{\cal C}|
\le |{\cal C}'| = \ip_{f'}(G')$. \qed

Consequently, we have the following result for $f$-isometric-path
numbers of connected block graphs.

\begin{thm}\label{thm:3} 
If $G$ is a connected block graph with a regular vertex labeling
$f$, then $\ip_f(G)=\lceil{s(G)}/{2}\rceil$, where
$s(G)=\sum_{v\in V(G)} f(v)$.
\end{thm}
{\bf Proof.}
The theorem is obvious when $G$ has only one vertex.
For the case when $G$ has more than one vertex,
we apply Theorem \ref{thm:2} repeatedly
until the graph becomes trivial.
Notice that the $s(G')=s(G)$ when apply Theorem \ref{thm:2}.
\qed

For the isometric-path-cover problem, we have

\begin{cor} \label{col:4} 
If $G$ is a connected block graph,
then $\ip(G)= \lceil{\nc(G)}/{2}\rceil$,
where $\nc(G)$ is the number of non-cut-vertices of $G$.
\end{cor}
{\bf Proof.}
The corollary follows from Theorem \ref{thm:3}
and the fact that $\ip(G)=\ip_f(G)$ for the regular vertex labeling $f$
with $f(v)=1$ if $v$ is a non-cut-vertex and $f(v)=0$ otherwise.
\qed

\section{Algorithm}

Based on Theorem \ref{thm:2}, we are able to design an algorithm
for the isometric-path-cover problem in block graphs. Notice that
we may only consider connected block graphs with regular vertex
labelings. To speed up the algorithm, we may modify Theorem
\ref{thm:2} a little bit so that each time a non-cut-vertex is
handled.

\begin{thm} \label{thm:5} 
Suppose $G$ is a block graph with a regular labeling $f$,
and $x$ is a non-cut-vertex in a block $B$
with exactly one cut-vertex $y$ or with no
cut-vertex  in which let $y$ be any vertex in $B-\{x\}$.
Let $G'=G-x$ with a regular vertex labeling $f'$
which is the same as $f$ except $f'(y)=f(y)+f(x)$.
Then $\ip_f(G)=\ip_{f'}(G')$.
\end{thm}
{\bf Proof.} The theorem follows from repeatedly applying Theorem
\ref{thm:2}. \qed

Now, we are ready to give the algorithm.

\bigskip

\noindent
{\bf Algorithm PG} Find the $f$-isometric-path number $\ip_f(G)$ of
                   a connected block graph.

\noindent
{\bf Input.} A connected block graph $G$ and a regular vertex labeling $f$.

\noindent
{\bf Output.} An optimal $f$-isometric-path cover $\cal C$ of $G$
              and $\ip_f(G)$.

\noindent
{\bf Method.}

\baselineskip 20pt

\begin{tabbing}
~1. ~~ \=  construct a stack $S$ which is empty at the beginning;\\
~2.    \>  let $G'\leftarrow G$;\\
~3.    \> {\bf while} ($G'$ has more than one vertex) {\bf do}\\
~4.    \> \quax \= choose a block $B$ with exactly one cut-vertex $y$ or with\\
       \>       \> \quax \= no cut-vertex in which case choose any $y\in B$;\\
~5.    \>       \> {\bf for} (all vertices $x$ in $B-\{y\}$) {\bf do}\\
~6.    \>       \> \quax \= $f(y)\leftarrow f(y)+f(x)$;\\
~7.    \>       \>       \> push $(x,y,f(x))$ into $S$;\\
~8.    \>       \>       \> $G' \leftarrow G'-x$;\\
~9.    \>       \> {\bf end for};\\
10.    \> {\bf end while};\\
11.    \> $\ip_f(G)\leftarrow \lceil{f(r)}/{2}\rceil$,
          where $r$ is the only vertex of $G'$;\\
12.    \> let $\cal C$ be the family of isometric paths
          containing $\ip(G)$ copies of the path $r$;\\
13.    \> {\bf while} ($S$ is not empty) {\bf do}\\
14.    \> \quax \= pop $(x,y,i)$ from $S$;\\
15.    \>       \> choose $i$ copies of path $P$ in $\cal C$
                   using $y$ as an end vertex;\\
16.    \>       \> {\bf if}   ($P=yx$)  {\bf then}\\
17.    \>       \> \quax   \= replace the $i$ copies of $P$
                              by $i$ copies of $x$ in $\cal C$;\\
18.    \>       \> {\bf if}   ($P=yz$ for some vertex $z$
                              in the block of $G$ containing $x$)  {\bf then}\\
19.    \>       \> \quax   \= replace the $i$ copies of $P$
                              by $i$ copies of $xz$ in $\cal C$;\\
20.    \>       \> {\bf if}   ($P=yQ$ where $Q$ has no vertices in
                              the block of $G$ containing $x$)  {\bf then}\\
21.    \>       \> \quax   \= replace the $i$ copies of $P$
                              by the $i$ copies of $xyQ$ in $\cal C$;\\
22.    \> {\bf end while}.
\end{tabbing}

\baselineskip 24pt

Algorithm {\bf PG} can be implemented in time linear to the number
of vertices and edges.  Notice that we can use the depth-first
search to find all blocks and cut-vertices of a graph, see
\cite{clr}.

\bigskip \noindent {\bf Acknowledgement.}  The authors thank the
referee for many constructive suggestions.

\end{document}